\documentclass{article}
\usepackage{amsmath}
\usepackage{graphicx}
\usepackage{latexsym}
\usepackage{hyperref}

\usepackage{amsfonts}
\usepackage[all]{xy}
\newcommand{\R}{\mathbb{R}}

\def\p{\partial}
\def\s{\sigma}
\def\na{\nabla}
\def\n{\nabla}
\def\ov{\overline}

\def\O{\Omega}
\def\H{\mathbf{H}}
\newcommand\be{\begin{eqnarray}}
\newcommand\en{\end{eqnarray}}
\def\en{\end{eqnarray}}
\def\div{{\mathrm{div}}}

\def\ino{\int_{\Omega}}
\def\inpo{\int_{\partial\Omega}}
\newcommand{\no}{\nonumber}

\def\D{\Delta}

\def\o{\omega}
\def\({\left(}
\def\){\right)}
\newtheorem{theorem}{Theorem}[section]

\newtheorem{corollary}[theorem]{Corollary}

\newtheorem{remark}[theorem]{Remark}

\begin{document}
\setcounter{page}{1}
\title{Estimates for eigenvalues of the Neumann and Steklov problems}
\author{Feng Du$^{a,b}$, Jing Mao$^{b,c}$,
Qiaoling Wang$^{d}$, Changyu Xia$^{d}$, Yan Zhao$^{b}$}
\date{}
\protect\footnotetext{\!\!\!\!\!\!\!\!\!\!\!\!{ MSC 2020: 35P15;
53C40; 58C40.}
 }
\maketitle ~~~\\[-15mm]

\begin{center}
{\footnotesize  $a$.  School of Mathematics and Physics Science,\\
Jingchu University of Technology, Jingmen, 448000, China\\ $b$.
Faculty of Mathematics and Statistics,\\
 Key Laboratory of Applied
Mathematics of Hubei Province, \\
Hubei University, Wuhan 430062, China\\
$c$. Department of Mathematics, Instituto Superior T\'{e}cnico, University of Lisbon,\\
Av. Rovisco Pais, 1049-001 Lisbon, Portugal\\
$d$. Departamento de Matem\'atica, Universidade de Brasilia,
70910-900-Brasilia-DF,
Brazil\\
 Emails:  defengdu123@163.com (F. Du), jiner120@163.com (J. Mao), \\
 wang@mat.unb.br (Q.
Wang), xia@mat.unb.br (C. Xia), hdxbzy@163.com (Y. Zhao).  }
\end{center}


\begin{abstract}

We prove Li-Yau-Kr\"oger type  bounds for Neumann-type eigenvalues
of the poly-harmonic operator and of the biharmonic operator on
bounded domains in a Euclidean space. We also prove sharp estimates
for lower order eigenvalues of a biharmonic  Steklov problem  and of
the Laplacian, which directly implies two sharp Reilly-type
inequalities for the corresponding first nonzero eigenvalue.

 \end{abstract}

\markright{\sl\hfill  F. Du, J. Mao, Q. Wang, C. Xia, Y. Zhao
\hfill}

\section{Introduction}
\renewcommand{\thesection}{\arabic{section}}
\renewcommand{\theequation}{\thesection.\arabic{equation}}
\setcounter{equation}{0} \label{intro}

Throughout this paper, let $\Omega$ be a bounded  domain with smooth
boundary $\p\O$ in the Euclidean $n$-space $\R^n$. We consider the
Neumann eigenvalue problem of the Laplacian $\ov{\D}$
\begin{eqnarray}\label{a1}
\left\{\begin{array}{ccc} -\ov{\D} u=\mu u\,&&~\mbox{in} ~~ \O, \\[2mm]
\frac{\p u}{\p\nu}=0 &&~~\mbox{on}~~\partial \O,
\end{array}\right.
\end{eqnarray}
 where $\frac{\p}{\p\nu}$ is the outward normal derivative on the boundary $\p\Omega$ w.r.t the outward unit  normal vector $\nu$. The
 system (\ref{a1}) can be used to describe the vibration of
 membrane and is also called the \emph{free membrane problem}. It is well known that this
problem has discrete spectrum $\{\mu_i\}_{i=1}^{\infty}$ diverging
to infinity and satisfying
 \begin{eqnarray*}
0=\mu_{1}(\Omega)<\mu_{2}(\Omega)\leq\mu_{3}(\Omega)\leq\cdots\uparrow
+\infty.
 \end{eqnarray*}

In \cite{AB},  Ashbaugh and Benguria  conjectured that
\begin{eqnarray}\label{a2}
\sum_{i=1}^{n}\frac{1}{\mu_{i+1}(\O)}\geq\frac{n}{\mu_2(B_\Omega)},
~~\mathrm{with~equality~if~and~only~if}~\Omega~\mathrm{is~a~ball},
\end{eqnarray}
where $B_\O$ is the ball of same volume as $\O$, $\mu_i(\O)$ is the
$i$-th Neumann eigenvalue on $\Omega$, and $\mu_2(B_\Omega)$ is the
first nonzero Neumann eigenvalue on $B_\Omega$. In \cite{WX3}, Wang
and Xia proved that
\begin{eqnarray}\label{a3}
\sum_{i=1}^{n-1}\frac{1}{\mu_{i+1}(\O)}\geq\frac{n-1}{\mu_2(B_\Omega)},
~~\mathrm{with~equality~if~and~only~if}~\Omega~\mathrm{is~a~ball},
\end{eqnarray}
which supports the above conjecture of  Ashbaugh and Benguria.

On the other hand,  corresponding to the Li-Yau's classical result
for Dirichlet eigenvalues of the Laplacian \cite{LY}, Kr\"oger
\cite{Kr} obtained the following inequality for the sum of the
Neumann eigenvalues
\begin{eqnarray}\label{a4}
\sum_{j=1}^k
\mu_{j}(\O)\leq(2\pi)^{2}\frac{n}{n+2}k^{\frac{n+2}{n}}\(\frac{1}{\omega_n|\Omega|}\)^{\frac{2}{n}},~~k\geq
1,
\end{eqnarray}
and the upper bound estimate for the $(k+1)$-th Neumann eigenvalue
\begin{eqnarray}\label{a5}
\mu_{k+1}(\O)\leq
(2\pi)^{2}\(\frac{n+2}{4\omega_n|\Omega|}\)^{\frac{2}{n}}k^{\frac{2}{n}},~~k\geq
0,
\end{eqnarray}
where $\omega_n$ denotes the volume of the unit ball in $\R^n$ and
$|\O|$ represents the volume of $\O$.

In this paper, we consider the following eigenvalue problem of the
poly-harmonic operator
\begin{eqnarray}\label{a6}
\left\{\begin{array}{ccc} \(-\ov{\D}\)^{p} u=\Gamma u\,&&~\mbox{in} ~~ \O, \\[2mm]
\frac{\p
u}{\p\nu}=\ov{\D}^{m}u=\frac{\p\ov{\D}^{m}u}{\p\nu}=\cdots=\ov{\D}^{2m-1}u=\frac{\p\ov{\D}^{2m-1}u}{\p\nu}=0
&&~~\mbox{on}~~\partial \O,~~\mathrm{when}~~p=2m,
\\[2mm] \frac{\p u}{\p\nu}=0 &&~~\mbox{on}~~\partial
\O,~~\mathrm{when}~~p=1,
\\[2mm]
~~~~\frac{\p
u}{\p\nu}=\frac{\p\ov{\D}^{m-1}u}{\p\nu}=\ov{\D}^{m}u=\frac{\p\ov{\D}^{m}u}{\p\nu}=\cdots=\ov{\D}^{2m-2}u\\
=\frac{\p\ov{\D}^{2m-2}u}{\p\nu}=0 &&~~\mbox{on}~~\partial
\O,~~\mathrm{when}~~p=2m-1, m>1,
\end{array}\right.
\end{eqnarray}
where  $p, m\in\mathbb{N}$ with $\mathbb{N}$ the set of all positive
integers.

\begin{remark}
\rm{ For the eigenvalue problem (\ref{a6}), we prefer to give some
facts as follows: \\
(1) Clearly, if $p=1$, then (\ref{a6}) degenerates into the
classical Neumann eigenvalue problem of the Laplacian, i.e.
(\ref{a1}). Based on this fact, by the abuse of terminology, we will
call (\ref{a6}) a Neumann-type
eigenvalue problem of the poly-harmonic operator $\(-\ov{\D}\)^{p}$. \\
 (2) In fact, when one considers the eigenvalue problem (\ref{a6}),
 the regularity assumption of $\partial\Omega$ should be made such that the embedding $H^{p}(\Omega)\subset
 L^{2}(\Omega)$ is compact. Our assumption on smoothness of $\partial\Omega$ here is
 stronger enough to ensure the compactness of this embedding.
 Actually, for Neumann eigenvalue problem of the Laplacian, i.e. the $p=1$ case of (\ref{a6}),
 the Lipschitz continuity assumption for $\partial\Omega$ is enough such that the embedding $H^{1}(\Omega)\subset
 L^{2}(\Omega)$ is compact. To avoid regularity argument, which is
 not so necessary for our main results of this paper, we have
 assumed that the boundary $\partial\Omega$ is smooth. \\
(3) Let $H^p(\Omega)$ denote the Sobolev space of functions in
$L^2(\Omega)$ with derivatives up to order $p$ in $L^2(\Omega)$. For
any $v,w\in H^p(\Omega)$, one can define an inner product
$\langle\cdot,\cdot\rangle$ as follows:
\begin{eqnarray*}
\langle v,w\rangle=\left\{
\begin{array}{ccc}
\int_{\Omega}\left[(\overline{\Delta}^{m}v)\cdot\overline{(\overline{\Delta}^{m}w)}+v\overline{w}\right],\,&&~
\mathrm{if}~p=2m,~m\in\mathbb{N}, \\[2mm]
\int_{\Omega}\left[\overline{\nabla}(\overline{\Delta}^{m-1}v)\cdot\overline{\overline{\nabla}(\overline{\Delta}^{m-1}w)}+v\overline{w}\right],
&&~~~~~~ \mathrm{if}~p=2m-1,~m\in\mathbb{N},
\end{array}
\right.
\end{eqnarray*}
where $\overline{\n}$ is gradient operator on $\Omega$, and volume
elements in the above integrals have been dropped.\footnote{~For
convenience, in the sequel we will drop the integral measures for
all integrals except it is necessary.} Here
\begin{eqnarray*}
\overline{\Delta}^{m}v=\sum\limits_{i_{1},i_{2},\cdots,i_{m}=1}^{n}\frac{\partial^{2m}v}{\partial
x^{2}_{i_{1}}\partial x^{2}_{i_{2}}\cdots\partial x^{2}_{i_{m}}},
\end{eqnarray*}
and $\overline{\Delta}^{m-1}$ can be defined similarly. The weak
version of (\ref{a6}) is then the variational problem
\begin{eqnarray*}
\int_{\Omega}(\overline{\Delta}^{m}u)\cdot\overline{(\overline{\Delta}^{m}w)}=\Gamma\int_{\Omega}u\overline{w},
~~\forall w\in H^{p}(\Omega)\qquad (if~ p=2m,~m\in\mathbb{N})
\end{eqnarray*}
 or
 \begin{eqnarray*}
\int_{\Omega}\overline{\nabla}(\overline{\Delta}^{m-1}u)\cdot\overline{\overline{\nabla}(\overline{\Delta}^{m-1}w)}=\Gamma\int_{\Omega}u\overline{w},
~~\forall w\in H^{p}(\Omega)\qquad (if~ p=2m-1,~m\in\mathbb{N})
\end{eqnarray*}
in unknowns $u\in H^{p}(\Omega)$ and $\Gamma\in\mathbb{R}$.
 It is not hard to verify that under the boundary conditions
proposed in (\ref{a6}), the poly-harmonic operator
$\(-\ov{\D}\)^{p}$ is self-adjoint w.r.t. the inner product
$\langle\cdot,\cdot\rangle$ defined as above. Then the standard
theory of self-adjoint compact operators tells us that the spectrum
of the eigenvalue problem (\ref{a6}) is real and discrete consisting
in a non-decreasing sequence
$$0=\Gamma_{1}(\Omega)<\Gamma_{2}(\Omega)\leq \Gamma_{3}(\Omega)\leq\cdots\uparrow +\infty,$$
where each eigenvalue (i.e., element in the discrete spectrum) is repeated with its multiplicity. \\
 (4) By the min-max principle, together with the divergence theorem and the boundary conditions in (\ref{a6}),
 it is not hard to know that the Rayleigh-Ritz type characterization of the $k$-th nonzero eigenvalue
$\Gamma_{k}(\Omega)$ is given as follows:
 \\ When $p=2m$,  $m\in\mathbb{N}$,
\begin{eqnarray}\label{a7}
\Gamma_{k}(\Omega)=\inf\left\{\frac{\ino(\ov{\D}^m u)^2}{\ino u^2
}\Bigg{|} u,\cdots, \ov{\D}^{2m-1} u\in H^2(\Omega), u\neq0,
\int_{\Omega}u u_{j}=0, j=1,\cdots,k-1\right\};
\end{eqnarray}
When $p=2m-1$, $m\in\mathbb{N}$,
\begin{eqnarray}\label{a8}
\Gamma_{k}(\Omega)=\inf\left\{\frac{\ino|\overline{\n}(\ov{\D}^{m-1}
u)|^2}{\ino u^2 }\Bigg{|} u,\cdots, \ov{\D}^{2m-2} u\in H^2(\Omega),
u\neq0, \int_{\Omega}u u_{j}=0, j=1,\cdots,k-1\right\},
\end{eqnarray}
where $u_j$ is the eigenfunction of the eigenvalue
$\Gamma_{j}(\Omega)$. Besides, the eigenfunction $u_{1}$ of
$\Gamma_{1}(\Omega)=0$ should be nonzero constant function.
 \\
(5) For $p=2$, by the min-max principle, one can deduce that the
$k$-th nonzero eigenvalue of (\ref{a6}) is not larger than the
square of the $k$-th nonzero Neumann eigenvalue of the Laplacian,
and so by using Schwarz inequality and (\ref{a3}), we infer that
\begin{eqnarray*}
\sum_{i=1}^{n-1}\frac{1}{\Gamma_{i+1}(\Omega)}\geq\frac{n-1}{\Gamma_{2}(B_\Omega)}.
\end{eqnarray*}
 (6) We have noticed that Provenzano \cite{lp} discussed the
 Dirichlet and the Neumann eigenvalues of the poly-harmonic operator
 $\(-\ov{\D}\)^{p}$ in the space $H^{p}(\Omega)$, with $\Omega\subset\mathbb{R}^n$ a bounded domain
 such that the embedding $H^{p}(\Omega)\subset L^{2}(\Omega)$ is compact,
  and successfully
 gave a relation between Dirichlet eigenvalues and Neumann
 eigenvalues therein -- the $(k+p)$-th Neumann eigenvalue is \emph{strictly
 less} than the $k$-th Dirichlet eigenvalue for all
 $k,p\in\mathbb{N}$. However, the boundary conditions proposed
 therein are different from the ones we have used in (\ref{a6}). In
 fact, Provenzano's boundary conditions proposed in \cite{lp} are
\begin{eqnarray*}
u=\frac{\partial
u}{\partial\nu}=\cdots=\frac{\partial^{p-1}u}{\partial\nu^{p-1}}=0~~\mathrm{on}~\partial\Omega
\qquad\qquad  (\mathrm{Dirichlet~case})
\end{eqnarray*}
and\footnote{~The Neumann boundary condition (\ref{NC}) can be
simplified as $\frac{\partial
u}{\partial\nu}\Big{|}_{\partial\Omega}=0$ for $p=1$,
$\frac{\partial^{2}u}{\partial\nu^{2}}=\mathrm{div}_{\partial\Omega}\left(\overline{\nabla}^{2}u\cdot\nu\right)_{\partial\Omega}
+\frac{\overline{\Delta}u}{\partial\nu}=0$ on $\partial\Omega$ for
$p=2$, where $\mathrm{div}_{\partial\Omega}$ is the surface
divergence on $\partial\Omega$, $\overline{\nabla}^{2}u$ is the
Hessian of $u$, and
$(\overline{\nabla}^{2}u\cdot\nu)_{\partial\Omega}$ stands for the
projection of $\overline{\nabla}^{2}u\cdot\nu$ to the tangent bundle
of $\partial\Omega$. As also mentioned in \cite{lp}, generally it is
``\emph{a quite involved task}" to write (\ref{NC}) explicitly for
$p\geq3$.}
\begin{eqnarray} \label{NC}
\mathcal{N}_{1}u=\mathcal{N}_{2}u=\cdots=\mathcal{N}_{p}u~~\mathrm{on}~\partial\Omega
\qquad\qquad  (\mathrm{Neumann~case}),
\end{eqnarray}
which leads to a truth that our eigenvalue problem (\ref{a6}) of the
operator $\(-\ov{\D}\)^{p}$ is different from those two investigated
in \cite{lp}.
 }
\end{remark}

We can prove the following Kr\"oger-type estimates for the
Neumann-type eigenvalue problem (\ref{a6}).

\begin{theorem} \label{th1}
Let $\Omega$ be a bounded connected domain, with smooth boundary
$\partial\Omega$, in the Euclidean $n$-space $\R^n$ and let
$\Gamma_j(\O)$ be the $j$-th eigenvalue of the system (\ref{a6}).
Then we have
\begin{eqnarray}\label{a9}
\sum_{j=1}^k
\Gamma_j(\O)\leq(2\pi)^{2p}\frac{n}{n+2p}k^{\frac{n+2p}{n}}\(\frac{1}{\omega_n|\Omega|}\)^{\frac{2p}{n}},~~k\geq
1,
\end{eqnarray}
and
\begin{eqnarray}\label{a10}
\Gamma_{k+1}(\Omega)\leq
(2\pi)^{2p}\(\frac{n+2p}{2p\omega_n|\Omega|}\)^{\frac{2p}{n}}k^{\frac{2p}{n}},~~k\geq
0,
\end{eqnarray}
where, as before, $\omega_n$ denotes the volume of the unit ball in
$\R^n$ and $|\O|$ denotes the volume of $\O$.
\end{theorem}

\begin{remark}
\rm{(1) Clearly, when $p=1$, the eigenvalue problem (\ref{a6})
reduces to (\ref{a1}), and then our upper bound estimates
(\ref{a9}), (\ref{a10})  become Kr\"oger's inequalities (\ref{a4})
and (\ref{a5}), respectively.
\\
 (2) Laptev \cite[Sections 2 and 3]{L} showed that Kr\"oger's estimates (\ref{a4})
and (\ref{a5}) are corollaries of general (sharp) trace inequalities
for convex functions of (self-adjoint) operators. Besides, in
\cite[Theorem 3.2]{L}, Laptev gave a lower bound for the counting
function of the spectrum of the Friedrich's extension
$B_{\mathcal{N}}$ of the differential operator
$B(D)=A^{\ast}(D)A(D)$ provided $B_{\mathcal{N}}$ has discrete
spectrum. Here the differential operator $A(D)$ is defined by
\begin{eqnarray*}
A(D)u(x)=\sum\limits_{\beta\leq l}A_{\beta}D^{\beta}u(x),\qquad u\in
C^{\infty}(\overline{\Omega},\mathbb{C}^{m}),~m\in\mathbb{N},~l\in\mathbb{N},
\end{eqnarray*}
where $\Omega\subset\mathbb{R}^{n}$ is an open set with its closure
$\overline{\Omega}$, the coefficients $A_{\beta}$ are $m\times
m$-matrices independent of $x\in\Omega$. Letting $m=1$,
$B(\xi)=|\xi|^{2}$, $l\in\mathbb{N}$, and then the operator
$B_{\mathcal{N}}$ coincides with the operator of the Neumann
boundary problem for the poly-harmonic operator $(-\Delta)^{l}$.
Hence, \cite[Theorem 3.2 and Corollary 3.3]{L} imply a lower bound
for the counting function of Neumann eigenvalues of $(-\Delta)^{l}$
-- see (3.5) in \cite{L}. Now, we claim that (\ref{a10}) is
equivalent to (3.5) in \cite{L}. In fact, let $N(\Gamma)= k$, where
$$N(\Gamma)=\sum_{\Gamma_{i}(\Omega)\leq
\Gamma}1=\mathop{\mathrm{sup}}_{\Gamma_{i}(\Omega)\leq \Gamma}i$$ is
the counting function, and then $\Gamma_{k+1}(\Omega)\geq\Gamma$.
Thus, we infer from (\ref{a10}) that
\begin{eqnarray*}
\Gamma\leq
(2\pi)^{2p}\(\frac{n+2p}{2p\omega_n|\Omega|}\)^{\frac{2p}{n}}k^{\frac{2p}{n}},\quad~~k\geq
0,
\end{eqnarray*}
which implies that
\begin{eqnarray*}
N(\Gamma)\geq
\frac{2p}{n+2p}\frac{1}{2\pi}\omega_n|\O|\Gamma^{\frac{n}{2p}}.
\end{eqnarray*}
}
\end{remark}

 Let $\ov{\D}$ and $\D$ be the
Laplace-Beltrami operators on $\O$ and $\p \O$, respectively. Let
$\ov{\n}$ and $\n$ be the gradient operators on $\O$ and $\p \O$
separately. Consider the following Neumann-type eigenvalue problem
of the biharmonic operator
\begin{eqnarray}\label{a11}
\left\{\begin{array}{ccc} \ov{\D}^2 u-\tau \ov{\D} u=\Lambda u\,&&~\mbox{in} ~~ \O, \\[2mm]
(1-\s)\frac{\p^2u}{\p\nu^2}+\sigma \ov{\D} u=0 &&~~\mbox{on}~~\partial \O,\\[2mm]
\tau\frac{\p u}{\p \nu}-(1-\s)\div_{\p \O}\(\ov{\nabla}^2
u\cdot\nu\)_{\partial\Omega}-\frac{\p\ov{\D} u}{\p\nu}=0
&&~~\mbox{on}~~\partial \O,
\end{array}\right.
\end{eqnarray}
where $\tau\geq 0$ and $\s\in(-1/(n-1),1)$ are two constants, and,
as before, $\div_{\p \Omega}$ denotes the tangential divergence
operator on $\p \Omega$, $\ov{\nabla}^2 u$ is the Hessian of $u$,
$(\overline{\nabla}^{2}u\cdot\nu)_{\partial\Omega}$ stands for the
projection of $\overline{\nabla}^{2}u\cdot\nu$ to the tangent bundle
of $\partial\Omega$. In this setting, the problem (\ref{a11}) has
discrete spectrum and all eigenvalues in the  spectrum can be listed
non-decreasingly as follows (see, e.g., \cite[Proposition 4.1]{C2})
\begin{eqnarray*}
0=\Lambda_{1}(\Omega)\leq\Lambda_{2}(\Omega)\leq\Lambda_{3}(\Omega)\leq\cdots\leq\cdots\uparrow
+\infty.
\end{eqnarray*}
 This problem is called the \emph{eigenvalue problem of free plate
under tension and with nonzero Poisson's ratio}, which for $n=2$ can
be used to describe the deformation of a planar material under
compression,
  $\tau,\s$ denote a parameter related to the tension and  a Poisson's ratio of the material, respectively.
 By the Rayleigh-Ritz characterization, the Neumann-type eigenvalues (if exist and with the abuse of terminology) of (\ref{a11})
  are given by (see, e.g., \cite{C2,LM} while \cite{BCL,C1} for the case $\sigma=0$)
\begin{eqnarray}\label{a12}
\Lambda_{k}(\Omega)=\mathop{\inf}_{0\neq u\in
H^2(\Omega)}\Bigg\{\frac{\ino\left[(1-\sigma)|\ov{\na}^2
u|^2+\s(\ov{\D}u)^2+\tau|\ov{\na} u|^2\right]}{\ino u^2}
\Bigg{|}\int_{\Omega}u u_j =0, j=1,\cdots,k-1 \Bigg\},
\end{eqnarray}
where
 $u_j$ is
an eigenfunction corresponding to the eigenvalue
$\Lambda_{j}(\Omega)$, and $|\ov{\na}^2
u|^{2}=\sum\limits_{i,j=1}^{n}\left(\frac{\partial^{2}u}{\partial
x_{i}\partial x_{j}}\right)^2$.

\begin{remark}
\rm{ (1) In \cite{C1,C2,LM}, the authors therein used the operator
$\mathrm{Proj}_{\partial\Omega}\left[(\overline{\nabla}^{2}u)\nu\right]$
to denote the projection  of $(\overline{\nabla}^{2}u)\nu$ onto the
space tangent to $\partial\Omega$, which obviously has the same
meaning as
$(\overline{\nabla}^{2}u\cdot\nu)_{\partial\Omega}$ here. \\
(2)  As before, let $B_\O\subset\R^n$ be the ball of same volume as
$\O$. When $\tau>0$, $\s=0$, Chasman \cite{C1} proved the following
isoperimetric inequality
$$\Lambda_1(\Omega)\leq\Lambda_1(B_\O), \quad~~\mathrm{with~equality~if~and~only~if}~\Omega~\mathrm{is~a~ball}.$$
When $\tau>0$, $\s\in (-1/(n-1), 1)$, Chasman \cite{C2} conjectured
that the above isoperimetric inequality is still true and
successfully proved a weaker version of it. Moreover, she also gave
numerical and analytic evidences to support this conjecture -- see
\cite[ Section 8]{C2} for details.}
\end{remark}

 When $\tau\geq 0$,
$0\leq\s<1$, for the eigenvalue problem (\ref{a11}), we can obtain
the following:
\begin{theorem} \label{th2}
Let $\Omega$, $|\O|$ and $\omega_n$ be as in Theorem \ref{th1}, and
let $\Lambda_j(\O)$ be the $j$-th eigenvalue of the system
(\ref{a11}).

(i) When $\tau\geq 0$ and $0\leq\s<1$, we have
\begin{eqnarray}\label{a13}
\sum_{j=1}^k
\Lambda_j(\O)\leq(2\pi)^{4}\frac{n}{(n+4)}k^{\frac{n+4}{n}}\(\frac{1}{\omega_n|\Omega|}\)^{\frac{4}{n}}
+\tau(2\pi)^{2}\frac{n}{(n+2)}k^{\frac{n+2}{n}}\(\frac{1}{\omega_n|\Omega|}\)^{\frac{2}{n}},\quad
~~k\geq 1;
\end{eqnarray}

(ii) When $\tau=0$ and $0\leq\s<1$, it holds
\begin{eqnarray}\label{a14}
\Lambda_{k+1}(\Omega)\leq
(2\pi)^{4}\(\frac{n+4}{4\omega_n|\Omega|}\)^{\frac{4}{n}}k^{\frac{4}{n}},\qquad~~k\geq
0;
\end{eqnarray}

(iii) When $\tau>0$ and $0\leq\s<1$, we have
\begin{eqnarray}\label{a15}
\Lambda_{k+1}(\Omega)\leq
\mathop{\mathrm{min}}_{r>2\pi\(\frac{k}{\omega_n|\O|}\)^{\frac{1}{n}}}
\frac{n\omega_n|\Omega|\(\frac{r^{n+4}}{n+4}+\tau\frac{r^{n+2}}{n+2}\)}{\omega_n|\Omega|
r^n-k(2\pi)^n},\qquad~~k\geq 0.
\end{eqnarray}
\end{theorem}

\begin{remark}
\rm{ (1) Recently, when $\tau\geq 0$, $\s=0$, Brandolini, Chiacchio
and Langford \cite{BCL} have already obtained upper bounds for the
sum of the first $k$ eigenvalues $\Lambda_{i}(\Omega)$ and for the
$(k+1)$-th eigenvalue $\Lambda_{k+1}(\Omega)$. Inspired by this fact
and our Theorem \ref{th2} here, together with the coercivity
argument for the sesquilinear form shown in \cite[Section 4]{C2},
the corresponding author, Prof. J. Mao, and his another collaborator
can also get the estimates (\ref{a13})-(\ref{a15}) under a more
general setting that $\tau\geq0$, $\sigma\in(-1/(n-1),1)$ -- see
\cite[Theorem 1.1 and Corollary 1.2]{LM} for details. Although
\cite{LM} has been published formally very recently, we still prefer
to remain Theorem \ref{th2} to emphasize and embody the origin and
continuity of our thought. \\
 (2) Clearly, if $\tau=0$ and $\sigma=1$, then
 \begin{eqnarray} \label{add2}
\left\{\begin{array}{ccc} \ov{\D}^2 u=\Lambda u\,&&~\mbox{in} ~~ \O, \\[2mm]
\ov{\D} u=\frac{\p\ov{\D} u}{\p\nu}=0 &&~~\mbox{on}~~\partial \O,
\end{array}\right.
\end{eqnarray}
which corresponds to the $p=2$ case of the eigenvalue problem
(\ref{a6}) except the boundary condition $\frac{\partial
u}{\partial\nu}\Big{|}_{\partial\Omega}=0$ missed. At end of
\cite[Section 4]{C2}, Chasman showed that for the eigenvalue problem
(\ref{add2}), all $H^{2}(\Omega)$ harmonic functions are
eigenfunctions with eigenvalue zero, and one has at least an
eigenvalue of infinite multiplicity. Based on this fact, we need to
expel $\tau=0$, $\sigma=1$ in Theorem \ref{th2} here and add the
boundary condition $\frac{\partial
u}{\partial\nu}\Big{|}_{\partial\Omega}=0$ in the previous
eigenvalue problem (\ref{a6}).
 }
\end{remark}

We also consider the following Steklov-type eigenvalue problem of
the biharmonic operator
\begin{eqnarray}\label{a16}
\left\{\begin{array}{ccc} \ov{\D}^2 u-\tau \ov{\D} u=0\,&&~\mbox{in} ~~ \O, \\[2mm]
(1-\s)\frac{\p^2u}{\p\nu^2}+\sigma \overline{\D} u=0 &&~~\mbox{on}~~\partial \O,\\[2mm]
\tau\frac{\p u}{\p \nu}-(1-\s)\div_{\p \Omega}\(\ov{\nabla}^2
u\cdot\nu\)_{\partial\Omega}-\frac{\p\ov{\D} u}{\p\nu}=\lambda u
&&~~\mbox{on}~~\partial \O,
\end{array}\right.
\end{eqnarray}
where $\tau,\sigma\in\mathbb{R}$ and other same symbols have the
same meanings as those in (\ref{a11}).

\begin{remark}
\rm{ (1) Li and Mao \cite[Theorem 2.1]{LM1} showed clearly that if
$\tau>0$ and $\sigma\in(-1/(n-1),1)$, the eigenvalue problem
(\ref{a16}) has the discrete spectrum and its elements (i.e.,
eigenvalues) can be listed non-decreasingly as follows
 \begin{eqnarray*}
0=\lambda_{0}(\Omega)<\lambda_{1}(\Omega)\leq\lambda_{2}(\Omega)\leq\cdots\leq\lambda_{k}(\Omega)\leq\cdots\uparrow\infty.
 \end{eqnarray*}
By means of variational principle, the Rayleigh-Ritz type
characterization of the $(k+1)$-th eigenvalue
$\lambda_{k+1}(\Omega)$ is given by
\begin{eqnarray}\label{a17}
\lambda_{k+1}(\Omega)=\mathop{\inf}_{0\neq u\in
H^2(\Omega)}\Bigg{\{}\frac{\int_{\Omega}\left[(1-\sigma)|\ov{\na}^2
u|^2+\s(\ov{\D}u)^2+\tau|\ov{\na} u|^2\right]}{\int_{\partial\Omega}
u^2
}\Bigg{|}\int_{\partial\Omega} u u_j=0, \nonumber\\
j=0,1,\cdots,k\Bigg{\}}, \qquad\qquad \qquad
\end{eqnarray}
where  $u_j$ is an eigenfunction corresponding to the eigenvalue
$\lambda_{j}(\Omega)$. Besides, the eigenfunction $u_{0}$ of
$\lambda_{0}(\Omega)=0$ should be nonzero constant function. \\
(2) When $\tau>0$, $\s=0$,
 Buoso and Provenzano \cite{BP} proved an isoperimetric inequality for the fundamental tone $\lambda_{1}(\Omega)$ of the system (\ref{a16})
 which states that
 $$\lambda_1(\Omega)\leq\lambda_1(B_\O),$$
with equality if and only if $\Omega$ is a ball. Here, as before,
 $B_\O\subset\R^n$ is the ball of same volume as $\O$. Very
recently,  Li and Mao \cite[Theorem 1.1]{LM1} showed that the above
isoperimetric inequality is still true for $\tau>0$ and
$\sigma\in(-1/(n-1),1)$, and moreover, the inequality can be
achieved when $\Omega$ is the ball $B_{\Omega}$. \\
 (3) For some other
estimates for $\lambda_{i}$'s, see, e.g., \cite{B,BC,BP,DMWX,XW}.
 }
\end{remark}

 Our next result is a sharp lower bound for the
sum of the reciprocals of the first $n$ nonzero eigenvalues of the
problem (\ref{a16}).
\begin{theorem} \label{th3}
Let $\Omega$ and $|\O|$ be as in Theorem \ref{th1}, and let
$\lambda_j(\O)$ be the $j$-th eigenvalue of the system (\ref{a16}).
When $\tau> 0$ and $\sigma\in(-1/(n-1),1)$, we have
\begin{eqnarray}\label{a18}
 \sum_{j=1}^n\frac{1}{\lambda_{j}(\Omega)} \geq \frac{|\p\O|^2}{\tau |\O|\inpo|\H|^2},
\end{eqnarray}
where $\H$ is the mean curvature vector of $\p\O$
 in $\R^n$,   $|\p\O|$ denotes the area of $\p\O$. Equality in (\ref{a18}) holds if and only
if $\O$ is a ball.
\end{theorem}

 Using the monotonicity of eigenvalues $\lambda_{i}$'s and
Theorem \ref{th3} immediately, we get
\begin{eqnarray*}
\frac{n}{\lambda_{1}(\Omega)}\geq\sum_{j=1}^n\frac{1}{\lambda_{j}(\Omega)}
\geq \frac{|\p\O|^2}{\tau |\O|\inpo|\H|^2},
\end{eqnarray*}
which directly implies the following Reilly-type eigenvalue
estimate.

\begin{corollary} \label{coro1}
Under the assumptions in Theorem \ref{th3}, we have
\begin{eqnarray*}
\lambda_{1}(\Omega)\leq n\tau\frac{|\O|}{|\p\O|^2}\inpo|\H|^{2},
\end{eqnarray*}
with equality holds if and only if $\O$ is a ball.
\end{corollary}

\begin{remark}
\rm{ Clearly, when the Reilly-type eigenvalue estimate in Corollary
\ref{coro1} attains the equality case, one also has
$\lambda_{1}(\Omega)=\lambda_{2}(\Omega)=\cdots=\lambda_{n}(\Omega)$.
}
\end{remark}

Our final result is  a sharp lower bound for the sum of the
reciprocals of the first $n$ nonzero eigenvalues of the Laplacian on
a closed submanifold immersed in a Euclidean space. Namely, we have:
\begin{theorem} \label{th4}
Let $M$ be an $n$-dimensional compact submanifold without boundary
isometrically immersed in  $\R^N$ and let $\eta_{j}(M)$ be the
$j$-th
 nonzero closed eigenvalue of the Laplacian on $M$. We have
\begin{eqnarray}\label{a19}
 \sum_{j=1}^{n}\frac{1}{\eta_{j}(M)} \geq \frac{|M|}{\int_M |\overline{\H}|^2},
\end{eqnarray}
where $\overline{\H}$ is the mean curvature vector of $M$
 in $\R^N$.
Moreover, when $n=N-1$, equality holds in (\ref{a19}) if and only if
$M$ is  a hypersphere of $\R^N$, and when $n<N-1$, if the
 equality holds in (\ref{a19}), then  $M$ is a  minimal  submanifold of some hypersphere of $\R^N$.
\end{theorem}

Using the monotonicity of nonzero closed eigenvalues $\eta_{i}$'s of
the Laplacian and Theorem \ref{th4} immediately, we get
\begin{eqnarray*}
\frac{n}{\eta_{1}(M)}\geq\sum_{j=1}^{n}\frac{1}{\eta_{j}} \geq
\frac{|M|}{\int_M |\overline{\H}|^2},
\end{eqnarray*}
which directly implies the following Reilly's eigenvalue estimate
(i.e., the main result of the  influential paper \cite{RR}).

\begin{corollary} \label{coro2}
Under the assumptions in Theorem \ref{th4}, we have
\begin{eqnarray*}
\eta_{1}(M)\leq\frac{n}{|M|}\int_M |\overline{\H}|^{2},
\end{eqnarray*}
and moreover, the equality holds implies the rigidity described as
in Theorem \ref{th4}.
\end{corollary}

\begin{remark}
\rm{ (1) Clearly, when the Reilly-type eigenvalue estimate in
Corollary \ref{coro2} attains the equality case, one also has
$\eta_{1}(M)=\eta_{2}(M)=\cdots=\eta_{n}(M)$, and if furthermore
$n=N-1$, then $\eta_{n+2}(M)>\eta_{n+1}(M)=\eta_{i}(M)$ for
$i=1,2,\cdots,n$, since the multiplicity of the first nonzero closed
eigenvalue of the Laplacian on any $n$-sphere in $\mathbb{R}^{n+1}$
is $n+1$ and the corresponding eigenfunctions are the restrictions
(to $n$-sphere) of $n+1$ coordinate functions of $\mathbb{R}^{n+1}$
(see, e.g., \cite[Chapter 2]{CI} for this fact).
 \\
 (2) Except Reilly's estimate for the first nonzero eigenvalue of
 the Laplacian (see \cite{RR} or Corollary \ref{coro2} here) and
 our Reilly-type estimate for the first nonzero eigenvalue of
 (\ref{a16}) --
 the Steklov-type eigenvalue problem of the biharmonic operator (see Corollary
 \ref{coro1}), some interesting Reilly-type estimates for the first
 nonzero eigenvalue of other type have also been obtained. For
 instance, Ilias and Makhoul \cite{IM} have obtained the Reilly-type
 estimate for the first nonzero Steklov eigenvalue of the Laplacian
 on compact submanifolds (with boundary) isometrically immersed
 in a Euclidean space; Du and Mao \cite{DM} have obtained the Reilly-type
 estimate for the first nonzero closed eigenvalue of the nonlinear
 $p$-Laplacian ($1<p<+\infty$) on compact submanifolds (without boundary) isometrically
 immersed into
a Euclidean space, a unit sphere, or even a projective space.
  }
\end{remark}

For convenience, in the sequel, we prefer to simplify the notations
for four types eigenvalues discussed in this paper, that is, we
separately write $\Gamma_{i}(\Omega)$, $\Lambda_{i}(\Omega)$,
$\lambda_{i}(\Omega)$ and $\eta_{i}(M)$ as $\Gamma_{i}$,
$\Lambda_{i}$, $\lambda_{i}$ and $\eta_{i}$. We also make an
agreement that these notations would be written completely if
necessary.

This paper is organized as follows. In Section \ref{S3}, we will
prove Li-Yau-Kr\"{o}ger type estimates for lower-order eigenvalues
of the Neumann-type eigenvalue problem (\ref{a6}) of the
poly-harmonic operator and the Neumann-type eigenvalue problem
(\ref{a11}) of the biharmonic operator. Two sharp extrinsic lower
bounds for the sum of the reciprocals of the first $n$ nonzero
eigenvalues of the Steklov-type eigenvalue problem (\ref{a16}) and
for the sum of the reciprocals of the first $n$ nonzero closed
eigenvalues of the Laplacian will be separately proven in Section
\ref{S3}.

\section{Li-Yau-Kr\"{o}ger type estimates}
\renewcommand{\thesection}{\arabic{section}}
\renewcommand{\theequation}{\thesection.\arabic{equation}}
\setcounter{equation}{0} \label{S2}

In this section, inspired by \cite{BCL,Kr,LY}, using the method of
Fourier transformation, together with the Rayleigh-Ritz type
characterizations (\ref{a7}), (\ref{a8}) and (\ref{a12}), we can
separately give the proofs of two Li-Yau-Kr\"{o}ger type estimates
by appropriately constructing trial functions.

\vspace{2mm}

First, we have:

$\\$\textbf{\emph{Proof of Theorem \ref{th1}}}.
 Let $\{\phi_j\}_{j=1}^{\infty}$ be
a  set of orthonormal eigenfunctions of the system (\ref{a6}), that
is,
\begin{eqnarray*}
\left\{\begin{array}{ccc} (-\ov{\D})^{p} \phi_j=\Gamma_j \phi_j\,&&~\mbox{in} ~~ \O,
 \\[2mm]
\frac{\p
\phi_{j}}{\p\nu}=\ov{\D}^{m}\phi_j=\frac{\p\ov{\D}^{m}\phi_j}{\p\nu}=\cdots=\ov{\D}^{2m-1}\phi_j=\frac{\p\ov{\D}^{2m-1}\phi_j}{\p\nu}=0
&&~~\mbox{on}~~\partial \O,~~\mathrm{when}~~p=2m,\\[2mm]
\frac{\p \phi_{j}}{\p\nu}=0 &&~~\mbox{on}~~\partial
\O,~~\mathrm{when}~~p=1,
\\[2mm]
~~~~~\frac{\p
\phi_{j}}{\p\nu}=\frac{\p\ov{\D}^{m-1}\phi_j}{\p\nu}=\ov{\D}^{m}\phi_j=\frac{\p\ov{\D}^{m}\phi_j}{\p\nu}=\cdots=\ov{\D}^{2m-2}\phi_j\\
=\frac{\p\ov{\D}^{2m-2}\phi_j}{\p\nu}=0
 &&~~\mbox{on}~~\partial \O,~~\mathrm{when}~~p=2m-1, m>1, \\[2mm]
\ino \phi_j\phi_l=\delta_{jl}.
\end{array}\right.
\end{eqnarray*}
Define
\begin{eqnarray*}
\Phi(x,y)=\sum_{j=1}^k \phi_j(x)\phi_j(y), ~~x,y\in \Omega,
\end{eqnarray*}
and let
\begin{eqnarray*}
\hat{\Phi}(z,y)=\frac{1}{(2\pi)^{\frac{n}{2}}}\ino \Phi(x,y)
e^{ix\cdot z}dx
\end{eqnarray*}
be the Fourier transform of $\Phi$ in the variable $x$.

Since
\begin{eqnarray*}(2\pi)^{\frac{n}{2}}\hat{\Phi}(z,y)=\sum_{j=1}^k \phi_j(y)\ino \phi_j(x) e^{ix\cdot z}dx
\end{eqnarray*}
is the orthogonal projection of the function $h_z(x)=e^{ix\cdot z}$
onto the subspace of $L^2(\Omega)$ spanned by $
\phi_1,\cdots,\phi_k$,
$\rho(z,y)=h_z(y)-(2\pi)^{\frac{n}{2}}\hat{\Phi}(z,y)$ can be used
as a trial function for $\Gamma_{k+1}$. Thus, we have from
(\ref{a7}) and (\ref{a8}) that
\begin{eqnarray*}
\Gamma_{k+1}\ino |\rho(z,y)|^2
dydz\leq\left\{\begin{array}{ccc}\ino|\ov{\D}^m_y\rho(z,y)|^2dydz,\qquad~~~~{\rm
if}\ p=2m,
\\[2mm]\ino|\overline{\n}_y(\ov{\D}^{m-1}_y\rho(z,y))|^2dydz,\qquad {\rm if}\  p=2m-1.
\end{array}\right.
\end{eqnarray*}
Integrating both sides of the above inequality over $B_r=\{z\in
\R^n\big{|}|z|<r\}$  yields
\begin{eqnarray}\label{b1}
\Gamma_{k+1}\leq
\left\{\begin{array}{ccc}\underset{r}{\inf}\frac{\int_{B_r}\ino|\ov{\D}^m_y\rho(z,y)|^2dydz}{\int_{B_r}\ino
|\rho(z,y)|^2 dydz},\ \ \ \ {\rm when\ } p=2m,\\
\underset{r}{\inf}\frac{\int_{B_r}\ino|\overline{\n}_y(\ov{\D}^{m-1}_y\rho(z,y))|^2dydz}{\int_{B_r}\ino
|\rho(z,y)|^2 dydz}, \ \ {\rm when\ }p=2m-1,
\end{array}
\right.
\end{eqnarray}
where $r>2\pi\(\frac{k}{\omega_n|\O|}\)^{\frac{1}{n}}$.
 Noticing
$|h_z(y)|=1$ and $\hat{\Phi}(z,y)=\sum_{j=1}^k
\phi_j(y)\hat{\phi_j}(z)$, we have
\begin{eqnarray}\label{b2}
\nonumber\int_{B_r}\ino |\rho(z,y)|^2 dydz&=&\int_{B_r}\ino
\left|h_z(y)-(2\pi)^{\frac{n}{2}}\hat{\Phi}(z,y)\right|^2dydz
\\\nonumber&=&\|h_z(y)\|^2-2(2\pi)^{\frac{n}{2}}\mathrm{Re}\(\int_{B_r}\ino h_z(y)\overline{\hat{\Phi}(z,y)}dydz\)+(2\pi)^n\|\hat{\Phi}(z,y)\|^2
\\\nonumber&=&\omega_n|\Omega|r^n-2(2\pi)^{\frac{n}{2}}\mathrm{Re}\(\sum_{j=1}^k\int_{B_r}\ino e^{iy\cdot z}\phi_j(y)\overline{\hat{\phi_j}(z)}dydz\)
\\\nonumber&&+(2\pi)^n\sum_{j,l=1}^k\int_{B_r}\ino \phi_j(y)\phi_l(y)\hat{\phi_j}(z)\overline{\hat{\phi_l}(z)}dydz
\\&=&\omega_n|\Omega|r^n-(2\pi)^n\sum_{j=1}^k\int_{B_r}|\hat{\phi_j}(z)|^2dz,
\end{eqnarray}
where $\|f\|^2=\int_{B_r}\ino |f|^2dydz.$

From $h_z(y)_{y_p}=\(e^{iy \cdot z}\)_{y_p}=iz_pe^{iy \cdot
z}=iz_ph_z(y)$, we get
 $$\overline{\Delta}_yh_z(y)= \sum_{p=1}^n h_z(y)_{y_py_p}=-|z|^2h_z(y),$$
 which gives
$$|\ov{\D}^m_yh_z(y)|^2=|z|^{4m},$$ and $$|\overline{\n}_y\ov{\D}^{m-1}_yh_z(y)|^2=|z|^{4m-2}.$$
Using integration by parts, we infer from the boundary condition of
(\ref{a6}) for $p=2m$ that
  \begin{eqnarray*}\ino \ov{\D}^m_yh_z(y)\overline{\ov{\D}^m_y\hat{\Phi}(z,y)}dy&=&\ino \ov{\D}^{m-1}_yh_z(y)\overline{\ov{\D}^{m+1}_y\hat{\Phi}(z,y)}dy
  \\&&-\inpo\(\ov{\D}^{m-1}_yh_z(y)\frac{\p\overline{\ov{\D}^m_y\hat{\Phi}(z,y)}}{\p\nu}
  -\overline{\ov{\D}^m_y\hat{\Phi}(z,y)}\frac{\p\ov{\D}^{m-1}_yh_z(y)}{\p\nu}\)dy
  \\&=&\ino \ov{\D}^{m-1}_yh_z(y)\overline{\ov{\D}^{m+1}_y\hat{\Phi}(z,y)}dy
  \\&=&\ldots \ldots \ldots \ldots \ldots \ldots \ldots \ldots \ldots\\
  &=&\ino h_z(y)\overline{\ov{\D}^{2m}_y\hat{\Phi}(z,y)}dy\\&&-\inpo\(h_z(y)\frac{\p\overline{\ov{\D}^{2m-1}_y\hat{\Phi}(z,y)}}{\p\nu}
  -\overline{\ov{\D}^{2m-1}_y\hat{\Phi}(z,y)}\frac{\p h_z(y)}{\p\nu}\)dy
  \\&=&\ino h_z(y)\overline{\ov{\D}^{2m}_y\hat{\Phi}(z,y)}dy.
  \end{eqnarray*}
So, when $p=2m$, we have
\begin{eqnarray}\label{b3}
\nonumber&&\int_{B_r}\ino|\ov{\D}^m_y\rho(z,y)|^2dydz
\\\nonumber&=&\int_{B_r}\ino\|\ov{\D}^m_yh_z(y)-(2\pi)^{\frac{n}{2}}\ov{\D}^m_y\hat{\Phi}(z,y)\|^2dydz
\\\nonumber&=&\|\ov{\D}^m_yh_z(y)\|^2-2(2\pi)^{\frac{n}{2}}\mathrm{Re}\(\int_{B_r}\ino \ov{\D}^m_yh_z(y)\overline{\ov{\D}^m_y\hat{\Phi}(z,y)}dydz\)
+(2\pi)^n\|\ov{\D}^m_y\hat{\Phi}(z,y)\|^2\\\nonumber&=&\frac{nr^{n+4m}}{n+4m}\omega_n|\Omega|
-2(2\pi)^{\frac{n}{2}}\mathrm{Re}\(\int_{B_r}\ino
h_z(y)\overline{\ov{\Delta}_y^{2m}\hat{\Phi}(z,y)}dydz\)
\\\nonumber&&+(2\pi)^n\int_{B_r}\ino\(\sum_{l_1=1}^k\ov{\D}^m_y\phi_{l_1}(y)\hat{\phi}_{l_1}(z)\)
\(\sum_{l_2=1}^k
\ov{\D}^m_y\phi_{l_2}(y)\overline{\hat{\phi}_{l_2}(z)}\)dydz
\\\nonumber&=&\frac{nr^{n+4m}}{n+4m}\omega_n|\Omega|
-2(2\pi)^{\frac{n}{2}}\mathrm{Re}\(\sum_{j=1}^k\Gamma_j\int_{B_r}\ino
e^{iy\cdot z}\phi_j(y)\overline{\hat{\phi}_j(z)} dydz\)
\\\nonumber&&+(2\pi)^n\sum_{l_1,l_2=1}^k\int_{B_r}\ino\ov{\Delta}_y^{2m}\phi_{l_1}(y)\hat{\phi}_{l_1}(z)
\phi_{l_2}(y)\overline{\hat{\phi}_{l_2}(z)}dydz
\\\nonumber
&=&\frac{nr^{n+4m}}{n+4m}\omega_n|\Omega|
-2(2\pi)^{n}\sum_{j=1}^k\Gamma_j\int_{B_r}|\hat{\phi}_j(z)|^2dz
\\\nonumber&&+(2\pi)^n\sum_{l_1,l_2=1}^k\Gamma_{l_1}\int_{B_r}\ino\phi_{l_1}(y)\hat{\phi}_{l_1}(z)
\phi_{l_2}(y)\overline{\hat{\phi}_{l_2}(z)}dydz\\
&=&\frac{nr^{n+4m}}{n+4m}\omega_n|\Omega|
-(2\pi)^{n}\sum_{j=1}^k\Gamma_j\int_{B_r}|\hat{\phi}_j(z)|^2dz.
\end{eqnarray}
Similarly, by using integration by parts, we can infer from the
boundary condition of (\ref{a6}) for $p=2m-1$  that
\begin{eqnarray*}&&\ino \overline{\n}_y(\ov{\D}^{m-1}_yh_z(y))\cdot\overline{\n}_y(\overline{\ov{\D}^{m-1}_y\hat{\Phi}(z,y)})dy\\&=&-\ino \ov{\D}^{m-1}_yh_z(y)\overline{\ov{\D}^{m}_y\hat{\Phi}(z,y)}dy+\inpo\ov{\D}^{m-1}_yh_z(y)\frac{\p\overline{\ov{\D}^{m-1}_y\hat{\Phi}(z,y)}}{\p\nu}dy
  \\&=&-\ino \ov{\D}^{m-1}_yh_z(y)\overline{\ov{\D}^{m}_y\hat{\Phi}(z,y)}dy
  \\&=&\ldots \ldots \ldots \ldots \ldots \ldots \ldots \ldots \ldots\\
  &=&-\ino h_z(y)\overline{\ov{\D}^{2m-1}_y\hat{\Phi}(z,y)}dy\\&&-\inpo\(h_z(y)\frac{\p\overline{\ov{\D}^{2m-2}_y\hat{\Phi}(z,y)}}{\p\nu}
  -\overline{\ov{\D}^{2m-2}_y\hat{\Phi}(z,y)}\frac{\p h_z(y)}{\p\nu}\)dy
  \\&=&-\ino h_z(y)\overline{\ov{\D}^{2m-1}_y\hat{\Phi}(z,y)}dy.
  \end{eqnarray*}
Therefore,
\begin{eqnarray}\label{b4}
\nonumber&&\int_{B_r}\ino|\overline{\n}_y(\ov{\D}^{m-1}_y\rho(z,y))|^2dydz
\\\nonumber&=&\int_{B_r}\ino|\overline{\n}_y(\ov{\D}^{m-1}_yh_z(y))-(2\pi)^{\frac{n}{2}}\overline{\n}_y(\ov{\D}^{m-1}_y\hat{\Phi}(z,y))|^2dydz
\\\nonumber&=&\|\overline{\n}_y(\ov{\D}^{m-1}_yh_z(y))\|^2-2(2\pi)^{\frac{n}{2}}\mathrm{Re}\(\int_{B_r}\ino \overline{\n}_y(\ov{\D}^{m-1}_yh_z(y))
\cdot\overline{\overline{\n}_y(\ov{\D}^{m-1}_y\hat{\Phi}(z,y))}dydz\)\\\nonumber&&
+(2\pi)^n\|\overline{\n}_y(\ov{\D}^{m-1}_y\hat{\Phi}(z,y))\|^2\\\nonumber&=&\frac{nr^{n+4m-2}}{n+4m-2}\omega_n|\Omega|
+2(2\pi)^{\frac{n}{2}}\mathrm{Re}\(\int_{B_r}\ino
h_z(y)\overline{\ov{\Delta}_y^{2m-1}\hat{\Phi}(z,y)}dydz\)
\\\nonumber&&-(2\pi)^n\int_{B_r}\ino\(\sum_{l_1=1}^k\ov{\D}^m_y\phi_{l_1}(y)\hat{\phi}_{l_1}(z)\)
\(\sum_{l_2=1}^k
\ov{\D}^{m-1}_y\phi_{l_2}(y)\overline{\hat{\phi}_{l_2}(z)}\)dydz
\\\nonumber&=&\frac{nr^{n+4m-2}}{n+4m-2}\omega_n|\Omega|
-2(2\pi)^{\frac{n}{2}}\mathrm{Re}\(\sum_{j=1}^k\Gamma_j\int_{B_r}\ino
e^{iy\cdot z}\phi_j(y)\overline{\hat{\phi}_j(z)} dydz\)
\\\nonumber&&-(2\pi)^n\sum_{l_1,l_2=1}^k\int_{B_r}\ino\ov{\Delta}_y^{2m-1}\phi_{l_1}(y)\hat{\phi}_{l_1}(z)
\phi_{l_2}(y)\overline{\hat{\phi}_{l_2}(z)}dydz
\\\nonumber&=&\frac{nr^{n+4m-2}}{n+4m-2}\omega_n|\Omega|
-2(2\pi)^{n}\sum_{j=1}^k\Gamma_j\int_{B_r}|\hat{\phi}_j(z)|^2dz
\\\nonumber&&+(2\pi)^n\sum_{l_1,l_2=1}^k\Gamma_{l_1}\int_{B_r}\ino\phi_{l_1}(y)\hat{\phi}_{l_1}(z)
\phi_{l_2}(y)\overline{\hat{\phi}_{l_2}(z)}dydz\\&=&\frac{nr^{n+4m-2}}{n+4m-2}\omega_n|\Omega|
-(2\pi)^{n}\sum_{j=1}^k\Gamma_j\int_{B_r}|\hat{\phi}_j(z)|^2dz.
\end{eqnarray}
Substituting (\ref{b2}) and (\ref{b3}) or (\ref{b4}) into (\ref{b1})
yields
\begin{eqnarray}\label{b5}
\Gamma_{k+1}\leq\mathop{\inf}_{r}\Bigg\{\frac{\frac{nr^{n+2p}}{n+2p}\omega_n|\Omega|
-(2\pi)^{n}\sum_{j=1}^k\Gamma_j\int_{B_r}|\hat{\phi}_j(z)|^2dz}{\omega_n|\Omega|r^n-(2\pi)^n\sum_{j=1}^k\int_{B_r}|\hat{\phi_j}(z)|^2dz}\Bigg\},
\end{eqnarray}
where the infimum is taken over
$r>2\pi\(\frac{k}{\omega_n|\O|}\)^{\frac{1}{n}}$. By Plancherel's
Theorem,
\begin{eqnarray}\label{b6}
c_j=\int_{B_r}|\hat{\phi_j}(z)|^2dz\leq 1~~~  \mathrm{for}
~~j=1,\cdots,k.
\end{eqnarray}
Combining (\ref{b5}) and (\ref{b6}), one gets
\begin{eqnarray*}
\Gamma_{k+1}\(\omega_n|\Omega|r^n-(2\pi)^n\sum_{j=1}^kc_j\)\leq\frac{n}{n+2p}\omega_n|\Omega|r^{n+2p}
-(2\pi)^{n}\sum_{j=1}^k\Gamma_jc_j,
\end{eqnarray*}
that is,
\begin{eqnarray*}
\Gamma_{k+1}\omega_n|\Omega|r^n-\frac{n}{n+2p}\omega_n|\Omega|r^{n+2p}\leq
(2\pi)^n\Gamma_{k+1}\sum_{j=1}^{k}c_j-(2\pi)^{n}\sum_{j=1}^{k}\Gamma_jc_j\leq
(2\pi)^{n}\sum_{j=1}^k\(\Gamma_{k+1}-\Gamma_j\).
\end{eqnarray*}
Solving the above inequality for $\sum_{j=1}^k \Gamma_j$, we have
\begin{eqnarray*}
(2\pi)^{n}\sum_{j=1}^k
\Gamma_j\leq\frac{n}{n+2p}\omega_n|\Omega|r^{n+2p}+\(k(2\pi)^{n}-\omega_n|\Omega|r^n\)\Gamma_{k+1}.
\end{eqnarray*}
Since $r>2\pi\(\frac{k}{\omega_n|\O|}\)^{\frac{1}{n}}$, then
$k(2\pi)^{n}-\omega_n|\Omega|r^n<0$, and we infer from the above
inequality that
$$\sum_{j=1}^k \Gamma_j\leq\frac{n\omega_n|\Omega|r^{n+2p}}{(n+2p)(2\pi)^{n}}.$$
Letting $r\rightarrow2\pi\(\frac{k}{\omega_n|\O|}\)^{\frac{1}{n}}$,
(\ref{a9}) follows.

Combining (\ref{b5}) and (\ref{b6}), we have
\begin{eqnarray}\label{b7-1}
\Gamma_{k+1}\leq\frac{\frac{nr^{n+2p}}{n+2p}\omega_n|\Omega|
}{\omega_n|\Omega|r^n-k(2\pi)^n}=F(r),\qquad
~~r>2\pi\(\frac{k}{\omega_n|\O|}\)^{\frac{1}{n}}.
\end{eqnarray}
Solving $F'(r)=0$, we get
$$r=2\pi\(\frac{(n+2p)k}{2p\omega_n|\Omega|}\)^{\frac{1}{n}}.$$
Taking the above value of $r$ into (\ref{b7-1}), we have
(\ref{a10}).
 \hfill $\Box$

\vspace{5mm}

At the end of this section, we also have:

$\\$ \textbf{\emph{Proof of Theorem \ref{th2}}}. Let
$\{\psi_j\}_{j=1}^{\infty}$ be the set of  orthonormal
eigenfunctions of the system (\ref{a11}), that is,
 \begin{eqnarray*}
\left\{\begin{array}{ccc} \ov{\D}^2 \psi_j-\tau \ov{\D} \psi_j=\Lambda_j \psi_j&&~\mbox{in} ~~ \O, \\[2mm]
(1-\s)\frac{\p^2\psi_j}{\p\nu^2}+\sigma \D \psi_j=0 &&~~\mbox{on}~~\partial \O,\\[2mm]
\tau\frac{\p \psi_j}{\p \nu}-(1-\s)\mathrm{div}_{\p \O}\(\ov{\nabla}^2 \psi_{j}\cdot\nu\)_{\partial\Omega}-\frac{\p\ov{\D} \psi_j}{\p\nu}=0 &&~~\mbox{on}~~\partial \O,\\[2mm]
\ino \psi_j \psi_l=0.
\end{array}\right.
\end{eqnarray*}
As in the proof of Theorem \ref{th1},  we know that
\begin{eqnarray*}(2\pi)^{\frac{n}{2}}\hat{\Psi}(z,y)=\sum_{j=1}^k \psi_j(y)\ino \psi_j(x) e^{ix\cdot z}dx
\end{eqnarray*}
is the orthogonal projection of the function $h_z(x)=e^{ix\cdot z}$
onto the subspace of $L^2(\Omega)$ spanned by $
\psi_1,\cdots,\psi_k$. Thus, we can use
$\varphi(z,y)=h_z(y)-(2\pi)^{\frac{n}{2}}\hat{\Psi}(z,y)$ as trial
function for $\Lambda_{k+1}$ to obtain
\begin{eqnarray*}
\Lambda_{k+1}\ino |\varphi(z,y)|^2
dydz\leq\ino\left[(1-\sigma)|\ov{\na}^2_y
\varphi(z,y)|^2+\s\|\ov{\D}_y\varphi(z,y)\|^2+\tau|\ov{\na}_y
\varphi(z,y)|^2\right]dydz.
\end{eqnarray*}
Integrating both sides of the above inequality over $B_r=\{z\in
\R^n\big{|}|z|<r\}$ yields
\begin{eqnarray}\label{b8}
\Lambda_{k+1}\leq
\mathop{\inf}_{r}\Bigg\{\frac{\int_{B_r}\ino\left[(1-\sigma)|\ov{\na}^2_y
\varphi(z,y)|^2+\s\|\ov{\D}_y\varphi(z,y)\|^2+\tau|\ov{\na}_y
\varphi(z,y)|^2\right]dydz}{\int_{B_r}\ino |\varphi(z,y)|^2
dydz}\Bigg\},
\end{eqnarray}
where $r>2\pi\(\frac{k}{\omega_n|\O|}\)^{\frac{1}{n}}$. By a similar
calculation to (\ref{b2}), we have
\begin{eqnarray}\label{b9}
\int_{B_r}\ino |\varphi(z,y)|^2
dydz=\omega_n|\Omega|r^n-(2\pi)^n\sum_{j=1}^k\int_{B_r}|\hat{\psi_j}(z)|^2dz.
\end{eqnarray}
Let $$P=\int_{B_r}\ino\left[(1-\sigma)|\ov{\na}^2_y
\varphi(z,y)|^2+\s\|\ov{\D}_y\varphi(z,y)\|^2+\tau|\ov{\na}_y
\varphi(z,y)|^2\right]dydz=P_1+P_2+P_3,$$ where
\begin{eqnarray*}
P_1&=&\int_{B_r}\ino\((1-\sigma)|\ov{\na}^2_y h_z(y)|^2+\s|\ov{\D}_yh_z(y)|^2+\tau|\ov{\na}_y h_z(y)|^2\)dydz,\\
P_2&=&-2(2\pi)^{\frac{n}{2}}\mathrm{Re}\Bigg\{\int_{B_r}\ino\Big((1-\sigma)\ov{\na}^2_y
h_z(y)\cdot\overline{\ov{\na}^2_y\hat{\Psi}(z,y)}
+\s\ov{\D}_y h_z(y)\overline{\ov{\D}_y\hat{\Psi}(z,y)}\\&&+\tau\ov{\na}_y h_z(y)\cdot\overline{\ov{\na}_y\hat{\Psi}(z,y)}\Big)dydz\Bigg\},\\
P_3&=&\int_{B_r}\ino\((1-\sigma)|\ov{\na}^2_y
\hat{\Psi}(z,y)|^2+\s|\ov{\D}_y\hat{\Psi}(z,y)|^2+\tau|\ov{\na}_y
\hat{\Psi}(z,y)|^2\)dydz.
\end{eqnarray*}
Since $|h_z(y)_{y_p}|=|z_p|$ and $|h_z(y)_{y_py_q}|=|z_p||z_q|$,
then $|\ov{\D}_yh_z(y)|=|z|^{2}$, $|\overline{\n}_y h_z(y)|=|z|$ and
$$|\ov{\na}^2 h_z(y)|^2=\sum_{p,q=1}^n|h_z(y)_{y_py_q}|^2=\sum_{p,q=1}^n|z_p|^2|z_q|^2=|z|^{4}.$$
 So, we have
\begin{eqnarray}\label{b10}
P_1=n\o_n|\O|\(\frac{r^{n+4}}{n+4}+\tau\frac{r^{n+2}}{n+2}\).
\end{eqnarray}
Integrating by parts and noticing
$\hat{\Psi}(z,y)=\sum_{j=1}^k\psi_j(y)\widehat{\psi_j}(z)$, it
follows that
\begin{eqnarray}\label{b11}
\nonumber
P_2&=&-2(2\pi)^{\frac{n}{2}}\mathrm{Re}\Bigg\{\int_{B_r}\ino\Big((1-\sigma)
h_z(y)\overline{\ov{\D}^2_y\hat{\Psi}(z,y)} +\s
h_z(y)\overline{\ov{\D}^2_y\hat{\Psi}(z,y)}
\\\nonumber&&-\tau h_z(y)\overline{\ov{\D}_y\hat{\Psi}(z,y)}\Big)dydz\Bigg\}\\
&=&-2(2\pi)^n\sum_{j=1}^n\Lambda_j
\int_{B_r}|\widehat{\psi_j}(z)|^2dz
\end{eqnarray}
and
\begin{eqnarray}\label{b12}
\nonumber P_3&=&\int_{B_r}\ino\((1-\sigma)|\ov{\na}^2_y
\hat{\Psi}(z,y)|^2+\s|\ov{\D}_y\hat{\Psi}(z,y)|^2+\tau|\ov{\na}_y
\hat{\Psi}(z,y)|^2\)dydz\\\nonumber
&=&\int_{B_r}\ino\hat{\Psi}(z,y)\overline{\(\ov{\D}_y^2-\tau\ov{\D}_y\)\hat{\Psi}(z,y)}dydz\\&=&
(2\pi)^n\sum_{j=1}^k\Lambda_j \int_{B_r}|\widehat{\psi_j}(z)|^2dz.
\end{eqnarray}
Combining (\ref{b8})-(\ref{b12}), we have
\begin{eqnarray}\label{b13}
\Lambda_{k+1}\leq\mathop{\mathrm{inf}}_{r}\Bigg\{\frac{\omega_n|\Omega|\(\frac{r^{n+4}}{n+4}+\tau\frac{r^{n+2}}{n+2}\)
-(2\pi)^{n}\sum_{j=1}^k\Lambda_j\int_{B_r}|\hat{\psi}_j(z)|^2dz}{\omega_n|\Omega|r^n-(2\pi)^n\sum_{j=1}^k\int_{B_r}|\hat{\psi_j}(z)|^2dz}\Bigg\}.
\end{eqnarray}
Letting
\begin{eqnarray}\label{b14}
c_j=\int_{B_r}|\hat{\psi_j}(z)|^2dz\leq 1~~~  \mathrm{for}
~~j=1,\cdots,k,
\end{eqnarray}
we deduce from (\ref{b13}) that
\begin{eqnarray*}
\Lambda_{k+1}\(\omega_n|\Omega|r^n-(2\pi)^n\sum_{j=1}^kc_j\)\leq
n\omega_n|\Omega|\(\frac{r^{n+4}}{n+4}+\tau\frac{r^{n+2}}{n+2}\)
-(2\pi)^{n}\sum_{j=1}^k\Gamma_jc_j,
\end{eqnarray*}
which implies that
\begin{eqnarray*}
\Lambda_{k+1}\omega_n|\Omega|r^n-n\omega_n|\Omega|\(\frac{r^{n+4}}{n+4}+\tau\frac{r^{n+2}}{n+2}\)\leq
(2\pi)^{n}\sum_{j=1}^k\(\Gamma_{k+1}-\Gamma_j\).
\end{eqnarray*}
Hence,
\begin{eqnarray*}
(2\pi)^{n}\sum_{j=1}^k \Lambda_j\leq
n\omega_n|\Omega|\(\frac{r^{n+4}}{n+4}+\tau\frac{r^{n+2}}{n+2}\)+\(k(2\pi)^{n}-\omega_n|\Omega|r^n\)\Gamma_{k+1}.
\end{eqnarray*}
Since $r>2\pi\(\frac{k}{\omega_n|\O|}\)^{\frac{1}{n}}$, we infer
from the above inequality that
$$\sum_{j=1}^k \Lambda_j\leq\(\frac{r^{n+4}}{n+4}+\tau\frac{r^{n+2}}{n+2}\)\frac{n\omega_n|\Omega|}{(2\pi)^{n}}.$$
Letting $r\rightarrow2\pi\(\frac{k}{\omega_n|\O|}\)^{\frac{1}{n}}$,
one gets (\ref{a13}).

Combining (\ref{b13}) and (\ref{b14}), we have
\begin{eqnarray}\label{b7}
\Lambda_{k+1}\leq\frac{\(\frac{r^{n+4}}{n+4}+\tau\frac{r^{n+2}}{n+2}\)\omega_n|\Omega|
}{\omega_n|\Omega|r^n-k(2\pi)^n}, \quad~~\forall
r>2\pi\(\frac{k}{\omega_n|\O|}\)^{\frac{1}{n}}.
\end{eqnarray}
Consequently, we have
\begin{eqnarray*}
\Lambda_{k+1}(\Omega)\leq
\mathop{\mathrm{min}}_{r>2\pi\(\frac{k}{\omega_n|\O|}\)^{\frac{1}{n}}}
\frac{n\omega_n|\Omega|\(\frac{r^{n+4}}{n+4}+\tau\frac{r^{n+2}}{n+2}\)}{\omega_n|\Omega|
r^n-k(2\pi)^n},\qquad~~k\geq 0.
\end{eqnarray*}
For the case $\tau=0$, solving $F'(r)=0$ yields
$$r=2\pi\(\frac{(n+2p)k}{2p\omega_n|\Omega|}\)^{\frac{1}{n}}.$$
Taking the above value of $r$ into (\ref{b7}), we have (\ref{a14}).
  \hfill $\Box$

\section{Reilly type estimates}
\renewcommand{\thesection}{\arabic{section}}
\renewcommand{\theequation}{\thesection.\arabic{equation}}
\setcounter{equation}{0} \label{S3}

In the last section, by using the QR-factorization theorem and the
variational principle, we can give the proofs of two sharp extrinsic
lower bounds for the sum of the reciprocals of the first $n$ nonzero
eigenvalues (given in Theorems \ref{th3} and \ref{th4}) by
constructing appropriately trial functions. In fact, we have already
used the method of QR-factorization (together with other approaches)
to try to get estimates for the sum of the reciprocals of the first
$n$ nonzero eigenvalues of prescribed eigenvalue problems (see,
e.g., \cite{LMW}).

\vspace{2mm}

First, we have:

$\\$\textbf{\emph{Proof of Theorem \ref{th3}}}. Let $x_1, \cdots,
x_n$ be the coordinate functions in $\R^n$. Since $\Omega$ is a
bounded domain in $\R^n$, we can regard $\p \Omega$ as a closed
hypersurface of $\R^n$ without boundary.

Let $u_j$ be an eigenfunction corresponding to the eigenvalue
$\lambda_{j}$ such that  $\{u_j\}_ {j=0}^{\infty}$ is an orthonormal
basis of $L^2(\partial \Omega)$, that is,
\begin{eqnarray*}
\left\{\begin{array}{ccc} \overline{\D}^2 u_j-\tau\overline{\D} u_j
=0&&~\mbox{in}~~\Omega, \\[2mm]
(1-\s)\frac{\p^2u_j}{\p\nu^2}+\sigma \ov{\D} u_j=0 &&~~\mbox{on}~~\partial \Omega,\\[2mm]
\tau\frac{\p u_j}{\p \nu}-(1-\s)\div_{\p \Omega}\(\overline{\nabla}^2 u_{j}\cdot\nu\)_{\partial\Omega}-\frac{\p\overline{\D} u_j}{\p\nu}=-\lambda_{j} u_i&&~~\mbox{on}~~\partial \Omega,\\[2mm]
\inpo u_i u_j =\delta_{ij}.
\end{array}\right.
\end{eqnarray*}
Observe that $u_0= 1/\sqrt{|\p \Omega|}$ is a constant. By
translating the origin appropriately, we can assume that \be
\label{c0} \int_{\p\O} x_i=0, \ i=1,\cdots,n,
 \en
 that is,
$x_i\perp u_0$. Next, we will show that a suitable rotation of axes
can be made so as to insure that
 \be\label{c1}
 \int_{\p\O} x_j
u_i=0,
 \en for $j=2,3,\cdots, n$ and $i=1,\cdots,j-1$. To see this,
define an $n \times n$ matrix $Q=\(q_{ji}\),$ where
$q_{ji}=\int_{\partial\Omega} x_j u_i$, for $i,j=1,2,\cdots,n.$
Using the orthogonalization of Gram and Schmidt (i.e.,
QR-factorization theorem), we know that there exist an upper
triangle matrix $T=(T_{ji})$ and an orthogonal matrix $U=(a_{ji})$
such that $T=UQ$, i.e.,
\begin{eqnarray*}
T_{ji}=\sum_{\gamma=1}^n x_{jk}q_{ki}=\inpo \sum_{k=1}^n  a_{jk}x_k
u_i =0,\ \  1\leq i<j\leq n.
\end{eqnarray*}
Letting $y_j=\sum_{k=1}^n  a_{jk}x_k$, we get
\begin{eqnarray}\label{c2}
\inpo y_j  u_i =\inpo \sum_{k=1}^n a_{jk}x_k u_i =0,\ \  1\leq
i<j\leq n.
\end{eqnarray}
Since $U$ is an orthogonal matrix,  $y_1, y_2, \cdots, y_n$ are also
coordinate functions on $\R^n$. Therefore, denoting these coordinate
functions still by $x_1, x_2,\cdots, x_n$, one can get (\ref{c1}).
From (\ref{c0}) and (\ref{c1}), one sees that
$x_j\bot\{u_0,u_1,\cdots,u_{j-1}\}$ in $L^2(\p\O)$.

It follows from the variational characterization (\ref{a17}) that
\begin{eqnarray*}
\lambda_{j}\inpo  x_j^2\leq \int_{\O} \(|\ov{\na}^2x_j|^2
+\tau|\ov{\na} x_j|^2\) = \tau |\O|, \qquad \ j=1,...,n,
\end{eqnarray*}
which implies that
\begin{eqnarray*}
 \sum_{j=1}^n\frac{1}{\lambda_{j}} \tau |\O|\geq \sum_{j=1}^n\inpo  x_j^2=\inpo |x|^2.
\end{eqnarray*}
Multiplying both sides of the above inequality by  $\inpo|\H|^2$,
and using the Schwarz inequality, we  obtain
\begin{eqnarray}\label{c3}
 \sum_{j=1}^n\frac{1}{\lambda_{j}} \tau |\O|\inpo|\H|^2\geq \inpo |x|^2\inpo|\H|^2\geq\(\inpo\langle x,\H\rangle\)^2=|\p\O|^2,
\end{eqnarray}
which gives (\ref{a18}).

If equality holds in (\ref{a18}), then all the inequalities  in
(\ref{c3}) should be equalities, which implies that $x=\kappa \H$
holds on $\p\O$ for some constant $\kappa\neq 0$. Thus, for any
tangent vector field $V$ on $\p \O$, we have $V(|x|^2) = 2\langle V,
x\rangle = 0$ and so $|x|$ and $|\H|$ are constants on $\p\O$. Since
$\p\O$
 is a closed hypersurface of $\R^n$, we conclude that $\p \O$
 is a round sphere.
 This completes the proof of Theorem \ref{th3}.\hfill $\Box$

 \vspace{5mm}

At the end, we also have:

$\\$\textbf{\emph{Proof of Theorem \ref{th4}}}. Denote by $\Delta$
and $\na$ the Laplacian and the gradient operator on $M$,
respectively. Without loss of generality, we can assume that $M$
does not lie in a hyperplane of $\R^N$. Let $x=(x_1, \cdots, x_N)$
be the position vector of $M$ in
 $\R^N$, and let $u_j$ be the normalized eigenfunction corresponding to the $j$-th nonzero  eigenvalue
$\mu_{j}$ of the Laplacian of $M$. By a similar discussion as in the
proof of Theorem \ref{th3}, we can assume that
$x_j\bot\{u_0,u_1,\cdots,u_{j-1}\}$ in $L^2(M)$. Then one has
\begin{eqnarray*}
\eta_{j}\int_M  x_j^2\leq \int_{M} |\na x_j|^2,\ j=1,\cdots,N,
\end{eqnarray*}
which implies that
\begin{eqnarray*}
 \sum_{j=1}^N\frac{1}{\eta_{j}}\int_{M}
|\na x_j|^2  \geq \sum_{j=1}^N\int_M  x_j^2=\int_M |x|^2.
\end{eqnarray*}
Since \be\no |\na x_j|^2\leq 1, \qquad \ \sum_{j=1}^N |\na x_j|^2=n,
\en we have \be \label{c4} \sum_{j=1}^N \frac{1}{\eta_j}|\na
x_j|^2&\leq&\sum_{j=1}^{n} \frac{1}{\eta_j}|\na
x_j|^2+\frac{1}{\eta_{n+1}}\sum_{A=n+1}^N|\na x_A|^2
\\ \no &=& \sum_{j=1}^{n}
\frac{1}{\eta_j}|\na
x_j|^2+\frac{1}{\eta_{n+1}}\(n-\sum_{i=1}^{n}|\na x_j|^2\)
\\ \no &\leq & \sum_{j=1}^{n}
\frac{1}{\eta_j}|\na x_j|^2+\sum_{i=1}^{n}\frac{1}{\eta_i}(1-|\na
x_i|^2)
\\ \no &=& \sum_{j=1}^{n}\frac{1}{\eta_j},
\en
 which gives
 \be \label{c5}\sum_{j=1}^{n}\frac{1}{\eta_j}|M| \geq \int_M |x|^2. \en
Multiplying both sides of the above inequality by $\int_M |\H|^2$,
and using the Schwarz inequality, we have
\begin{eqnarray}\label{c6}
 \sum_{j=1}^n\frac{1}{\eta_{j}} |M|\int_M |\overline{\H}|^2\geq \int_M |x|^2\int_M|\overline{\H}|^2\geq\(\int_M \langle x,\overline{\H}\rangle\)^2=|M|^2,
\end{eqnarray}
which implies that (\ref{a19}) is true.

If equality holds in (\ref{a19}), then equalities hold in all of the
above inequalities, which implies that
\begin{eqnarray*}
\eta_{1}=\cdots=\eta_{N}\equiv C,
\end{eqnarray*}
\begin{eqnarray*}\label{c12}
\overline{\D} x_j =-C x_j,  j=1,\cdots, N, \ \ {\rm on}\ \ M,
\end{eqnarray*}
and $x=\kappa \overline{\H}$ hold on $M$ for some constant
$\kappa\neq 0$.  From these facts, we know that $|x|$ and
$|\overline{\H}|$ are constants on $M$. Therefore, when $n=N-1$, $M$
is a hypersphere, and when $n<N-1$, $M$ is a minimal  submanifold of
some hypersphere of $\R^N$.    \hfill $\Box$

\section*{Acknowledgments}
 F. Du was
supported by Hubei Key Laboratory of Applied Mathematics (Hubei
University), Research Team Project of Jingchu University of
Technology (Grant No. TD202006) and Research Project of Jingchu
University of Technology (Grant Nos. YB202010, ZX202002, ZX202006).
J. Mao was supported in part by the NSF of China (Grant No.
11801496), the Fok Ying-Tung Education Foundation (China) and Hubei
Key Laboratory of Applied Mathematics (Hubei University). He wants
to thank the Department of Mathematics, IST, University of Lisbon
for its hospitality during his visit from September 2018 to
September 2019. Q. Wang was supported by CNPq, Brazil (Grant No.
307089/2014-2). C. Xia was supported by CNPq, Brazil (Grant No.
306146/2014-2).

\end{document}